\documentstyle[12pt]{article}
\title{TOWARD A GENERAL THEORY OF TRANSMUTATION}
\author{Amin Boumenir\thanks{K.F.U.P.M., Dhahran, Saudi Arabia,
visiting University of Illinois}
\and Robert Carroll\thanks{University of Illinois, Urbana,
IL 61801}}
\date{December, 1994}

\begin{document}
\maketitle
\begin{abstract} A general construction of transmutation operators
is developed for selfadjoint operators in Gelfand triples.  Theorems
regarding analyticity of generalized eigenfunctions and Paley-Wiener
properties are proved.
\end{abstract}

\section{Introduction}
The idea of transmutation operator (or transformation
operator) B such that BP = QB for P and Q ordinary differential
operators goes back to Gelfand, Levitan, Marchenko, Naimark, et.
al. in the early 1950's (cf. [23;25;28]).  It was picked up again
by Delsarte and Lions, who established some fundamental ideas (cf.
[26;31]), and subsequently it was developed in many directions
(see e.g. [2;10-14;17;22;31]).  In this article we indicate some
constructions of a general nature which will be further enhanced
in subsequent papers.  We develop the theory via selfadjoint operators
in Gelfand triples and give some constructions of
transmutation operators with various domains.
Then various properties such as analyticity
of generalized eigenfunctions and Paley-Wiener properties are discussed,
with results of various kinds.

\section{Background}
A typical background situation involves $P = -D^2$ and
$Q = -D^2 + q$ (q real) on $[0,\infty),\,\,D \sim \partial_x$
where e.g. $q \in C^0[0,\infty)$ and $\int^{\infty}_{-\infty}
(1 + x^2)|q|dx < \infty$.  This is a typical inverse scattering
situation and we denote by $\phi = Coskx$ and $\psi$ the generalized
eiqenfunctions satisfying
$$ P\phi = k^2\phi;\, Q\psi = k^2\psi;\,\phi(0,k) = \psi(0,k) = 1;\,
\phi'(0,k) = \psi'(0,k) = 0
\eqno(2.1)$$
\medskip
\noindent Here we will also write $\lambda = k^2$ and $k = \sqrt{\lambda}$,
depending on context, with abuse of notation such as $\psi(x,k)
\sim \psi(x,\lambda)$ when no confusion can arise, and one notes
that $\phi \sim \phi^P_k$ with $\psi \sim \phi^Q_k$ in the notation
of [11-14].  More general initial conditions $h\psi(0,k) - \psi'(0,k)
= 0$ can also be envisioned.  Then (cf. [11;12;28] for details) one
can produce by PDE techniques or by Paley-Wiener theory a triangular
kernel K(x,y) such that
$$ \psi(x,\lambda) = \phi(x,\lambda) + \int^x_0K(x,t)\phi(t,\lambda)dt
= (B\phi(\cdot,\lambda))
\eqno(2.2)$$
\medskip
\noindent This can be written $\psi(x,\lambda) = (B\phi(\cdot,\lambda))
= <\beta(x,t),\phi(t,\lambda)>$ for $\beta(x,t) = \delta(x-t) +
K(x,t)$ and we will write, for suitable f,
$$ (B^{*}f)(t) = <\beta(x,t),f(x)> = f(t) + \int_t^{\infty}K(x,t)f(x)dx
\eqno(2.3)$$
\medskip
\noindent Further in the same way one can transmute in the opposite
direction via
$$ \phi(x,\lambda) = \psi(x,\lambda) + \int_0^xL(x,t)\psi(t,\lambda)dt
= ({\cal B}\psi(\cdot,\lambda))
\eqno(2.4)$$
\medskip
\noindent Thus we will write BP = QB and ${\cal B}Q = P{\cal B}$
with $\gamma(x,y) = \delta(x-y) + L(x,y)$.\\
\medskip
\indent We emphasize that $\phi,\,\,\psi\,\not\in L^2_x$ and we are
not at the moment dealing with an $L^2$ theory; the brackets $<\,,
\,>$ denote suitable distribution pairings.  Assume now that Q has only
continuous spectrum (P does of course) and consider transforms
(for suitable f)
$$ {\cal P}f(k) = {\cal C}f(k) = \int_0^{\infty}Coskx\,f(x)dx;
\eqno(2.5)$$
$${\cal P}^{-1}F(x) = \frac{2}{\pi}\int_0^{\infty}Coskx\,F(k)dk;
{\cal Q}f(k) = \int_0^{\infty}f(x)\psi(x,k)dx$$
\medskip
\indent Next it is shown in [11;12;27] for example that there is a
generalized spectral function $R^Q \in Z'$ and a Parseval formula
$$ <f,g> = <R^Q,{\cal Q}f{\cal Q}g>
\eqno(2.6$$
\medskip
\noindent for functions $f,g\,\in K^2 = \{f \in L^2(0,\infty)$ with
compact support\}.  Here ${\cal P}K^2 = {\cal C}K^2 = \cup
{\cal C}K^2(\sigma)$ for ${\cal C}K^2(\sigma) = \{$even entire
$\hat{f}(k) = {\cal P}f$ with $\hat{f} \in L^2$ for k real and
$|\hat{f}(k)| \leq cexp(\sigma|Im(k)|)$ via $f \in K^2(\sigma)$
or $supp(f)\subset [0,\sigma]\}$, and $Z = \cup Z(\sigma)$ for
$Z(\sigma) = \{$even entire functions g(k) with $g \in L^1$ for
k real and $|g(k)|\leq cexp(\sigma|Im(k)|)\}$.  Z has a countable
union topology as in [19] and $Z'$ is its dual.  From this one
obtains an inversion
$$ {\cal Q}^{-1}F(x) = <R^Q,F(k)\psi(x,k)>
\eqno(2.7)$$
\medskip
\noindent This all leads to factorizations
$$ {\cal P}B^{*} = {\cal Q};\,\,{\cal Q}{\cal B}^{*} = {\cal P}
\eqno(2.8)$$
\medskip
\noindent Further let us define operators, when they make
sense (here $R^P \sim
\frac{2}{\pi}dk$)
$$ \tilde{{\cal Q}}F(x) = <F(k)\psi(x,k),R^P>; \tilde{{\cal P}}
F(x) = <F(k)\phi(x,k),R^Q>
\eqno(2.9)$$
\medskip
\noindent
Then e.g. via formal representations
$$ \beta(y,x) = <\phi(x,k)\psi(y,k),R^P>; \gamma(x,y) = <\phi(x,k)
\psi(y,k),R^Q>
\eqno(2.10)$$
\medskip
\noindent one can write
$$ B = \tilde{{\cal Q}}{\cal P};\, {\cal B} = \tilde{{\cal P}}{\cal Q};\,
B^{*} = {\cal P}^{-1}{\cal Q};\, {\cal B}^{*} = {\cal Q}^{-1}{\cal P}
\eqno(2.11)$$
\medskip
\noindent Such formulas were applied to many operators P and Q in
[11;12;14], involving both singular and nonsingular situations, and
many explicit formulas for kernels etc. were obtained in terms of
special functions.  It was most often the case the the generalized
spectral functions $R^Q$ and $R^P$ were in fact measures $d\Gamma_Q$ and
$d\Gamma_P$ in which case one can rewrite the spectral pairings
etc. as integrals.  In particular (2.6) becomes
$$ \int fg\,dx = \int {\cal Q}f{\cal Q}g\,d\Gamma_Q
\eqno(2.12)$$
\medskip
\noindent Further the measures $d\Gamma_P,\,d\Gamma_Q$ were frequently
absolutely continuous with say $d\Gamma_P = \gamma_Pdk$ and this
will set the stage for our presentation in section 3.

\section{Transmutation for certain selfadjoint operators}

The theory described in section 2 was based on differential operators
$Q = -D^2 + q$ and the associated Paley-Wiener theory for example.
We want to deal now with a more general situation where less is
assumed a priori and which will include the type of situation
described in section 2.  Thus take a densely defined selfadjoint
operator Q in $L^2(dM_Q) = L^2_Q$ with a simple spectrum (cf. [1]).
This entails no basic loss in generality since for finite multiplicity
one could decompose Q as a direct sum of operators with simple
spectrum.  For $\lambda \in sp(Q) = \sigma_Q$ there exists a sequence
$f_n \in D(Q) \subset L^2_Q$ such that $\|f_n\| = 1$ and
$lim\,\|Qf_n - \lambda f_n\| = 0$.  If $f_n \to f$ in $L^2_Q$ then
$\lambda \in \sigma_Q^d$ (discrete spectrum) and $Qf = \lambda f$
with $f \in L^2_Q$.  If $f_n$ does not converge then $\lambda \in
\sigma_Q^c$ (continuous spectrum).  We imagine
then that $L^2_Q \hookrightarrow \Phi'$ with compact embedding
(e.g. think of Gelfand triples or rigged Hilbert spaces as in
[3;4;13;18-20;28]) and then $f_n \to F \in \Phi'$.  Thus we will
have a solution $F \sim \psi(x,\lambda) \in \Phi'$ of $Q\psi =
\lambda\psi$ for $\lambda \in \sigma_Q$.  As to constructing such
$\Phi$ one recalls (cf. [1]) that given a selfadjoint operator Q in a
Hilbert space H with simple spectrum there is a vector $h \in D(Q)
\subset H$ such that $Q^kh$ is defined for all k and the linear
envelope $\Phi = \{Q^kh\}$ is dense in H.  Put on $\Phi$ if possible
a topology such that $i: \,\Phi\hookrightarrow H$ is compact and
embed H in the antidual $\Phi'$ via
$\phi \to (\phi,h) = <L_h,\phi>$ (for convenience from now on
we will think of real Hilbert spaces without loss of generality
-cf. [12;14]).  Then $\Phi \subset D(Q)$ and $Q:\,\Phi \to \Phi$ is
continuous (variations on this are indicated below).  Note
$H_n = \{\sum_{|k|\leq n}a_kQ^k h\}$ is finite dimensional, hence
nuclear, and $\Phi = \overrightarrow{lim}\,H_n$ is nuclear
with $\Phi\hookrightarrow H$.  However $i:\,\Phi \to H$ is not
a priori compact or Hilbert-Schmidt without further hypotheses.  The
theory of rigged spaces (cf. [3;19;20;29]) then provides a measure
$d\Gamma_Q$ with ($f \in \Phi$)
$$ {\cal Q}f(\lambda) = \hat{f}(\lambda) = <f(x),\psi(x,\lambda);
\eqno(3.1)$$
$$ f(x) = \int_{-\infty}^{\infty}{\cal Q}f(\lambda)\psi(x,\lambda)
d\Gamma_Q(\lambda); \|f\|^2 = \int_{-\infty}^{\infty}|{\cal Q}
f(\lambda)|^2 d\Gamma_Q(\lambda)$$
\medskip
\noindent Thus ${\cal Q}:\,\Phi \to L^2(d\Gamma_Q)$ is continuous
and this can be extended to an isometry $\bar{{\cal Q}}:\,H \to
L^2(d\Gamma_Q)$,
which we usually denote again by ${\cal Q}$.  We observe now that
given the spectral measure $d\Gamma_Q$ and $\psi(x,\lambda) \in$
(some $\hat{\Phi}'$) with $H \sim L^2_Q$ but no $\hat{\Phi}$ in sight,
it would be natural to expect ${\cal S}_{\lambda} \subset
L^2(d\Gamma_Q)$ (${\cal S}$ = Schwartz space).  Then from the property
$Q\psi = \lambda \psi$ one has for $F \in {\cal S}$ and $f(x) =
{\cal Q}^{-1}F = \int F(\lambda)\psi(x,\lambda)d\Gamma_Q$,
$$ Q^n f = \int F\,Q^n\psi\,d\Gamma_Q = \int(\lambda^n F)\psi d\Gamma_Q
\eqno(3.2)$$
\medskip
\noindent which is well defined since $\lambda^n F \in {\cal S}$ again.
It is not clear however how $\tilde{\Phi} = {\cal Q}^{-1}{\cal S}$
is related to $\Phi$ for $\Phi$ constructed above or to a putative
$\hat{\Phi}$.  In any case $\tilde{\Phi}$ will be a perfect space
(bounded sets are relatively compact) with the topology defined via
seminorms
$$ \|\tilde{\phi}\|_p = |{\cal Q}\tilde{\phi}|_p =
sup_{q\leq p}|(1 + \lambda^{2p})\partial^q_{\lambda}
({\cal Q}\tilde{\phi})(\lambda)|
\eqno(3.3)$$
\medskip
\noindent (i.e. with the topology induced by ${\cal S}$ and ${\cal Q}$).
Further $\tilde{\Phi}\hookrightarrow H$ with compact embedding and
$\tilde{\Phi}\hookrightarrow H  \hookrightarrow \tilde{\Phi}'$ with
$Q\tilde{\Phi}\subset\tilde{\Phi}$.  Hence one knows there exists
a generalized eigenfunction $\tilde{\psi}(x,\lambda) \in \tilde{\Phi}'$
and a measure $d\tilde{\Gamma}_Q$ with formulas (3.1) for $\tilde
{{\cal Q}}f(\lambda) = <f(x),\tilde{\psi}(x,\lambda)>$.  One has
$\tilde{{\cal Q}}\tilde{\Phi} = {\cal S} \subset L^2(\tilde{\Gamma}_Q)$
but generally one may not have $\psi(x,\lambda) \in \tilde{\Phi}'$.
This is now circumvented by transferring the theory of Q to
$\tilde{\Phi} \hookrightarrow H \hookrightarrow \tilde{\Phi}'$ and
using $\tilde{\psi}(x,\lambda)$ in place of $\psi(x,\lambda)$.
Consequently, given a selfadjoint Q in H with simple spectrum, once we
have a generalized eigenfunction $\psi(x,\lambda)$ with formulas
(3.1) we can produce a Gelfand triple $\tilde{\Phi} \hookrightarrow H
\hookrightarrow \tilde{\Phi}'$ with an isomorphic theory.  This leads
us to work with the class of operators (s.a. $\sim$ selfadjoint)
$$ {\cal A} = \{s.a.\,\, Q\,\, in\,\,H \sim
L^2_Q\,\,with\,\,(3.1),
\,\,\Phi \hookrightarrow H \hookrightarrow \Phi',\,\,\Phi \subset
D(Q)\,\,dense,\,\,Q\Phi\subset\Phi\}
\eqno(3.4)$$
\medskip
\noindent
where $\Phi$ is to be dense in D(Q) with graph norm and in H.\\
\medskip
\indent Now given $Q_1,\,Q_2 \in {\cal A}$ (based on $\Phi_i
\hookrightarrow H_i \hookrightarrow \Phi'_i$) let us assume first
$\sigma_1 = \sigma_2$ for simplicity and write for $f \in \Phi_i$
$$ \hat{f}_i(\lambda) = <f(x),\psi_i(x,\lambda)>;\,
f(x) = \int\hat{f}_i(\lambda)\psi_i(x,\lambda)d\Gamma_i(\lambda)
\eqno(3.5)$$
\medskip
\noindent (recall we are using real $L^2$ spaces for convenience -
the corresponding results for complex spaces will follow as indicated
in [12;14]).  Denote by ${\cal Q}_i$ the maps indicated in (3.1) so
${\cal Q}_i$ extends to an isometry ${\cal Q}_i: H_i \to L^2(d\Gamma_i)\,\,
(H_i \sim L^2_{Q_i} = L^2(dM_i))$.  Define now a map
$$ {\bf Q}: L^2(d\Gamma_1) \to L^2(d\Gamma_2):\,\hat{f}_1(\lambda)
\to {\bf Q}\hat{f}_1(\lambda) = \hat{f}_2(\lambda)
\eqno(3.6)$$
\medskip
\noindent Here $f \in H_1\cap H_2$ and $D({\bf Q}) = \{F \in
L^2(d\Gamma_1);\,\,{\cal Q}_1^{-1}F \in H_1\cap H_2\}$ with
\begin{eqnarray}
\setlength{\jot}{4.5pt}
\hat{f}_1 \hspace{4pt} & \stackrel{{\bf Q}}{\rightarrow} \hspace{4pt} &
{\bf Q}\hat{f}_1 = \hat{f}_2\nonumber\\
{\cal Q}_1\uparrow \hspace{4pt} & {\cal Q}_2\nearrow \hspace{4pt} &
\Uparrow {\cal Q}_1\nonumber\\
f\in H_1\cap H_2 \hspace{4pt} & \stackrel{V}{\Rightarrow} \hspace{4pt} &
Vf \in H_1\nonumber
\end{eqnarray}
\medskip
\noindent We will assume here that $H_1\cap H_2$ is dense in $H_i$.
The double arrows indicate a formal relation to the machinery of
section 1 where $V = {\cal Q}_1^{-1}{\cal Q}_2 \sim B^{*}$ in (2.10).
Here we think of transmutations $B:\,Q_1 \to Q_2,\,\,BQ_1 = Q_2B$,
acting on $f \in D(Q_1)$ with $Bf \in D(Q_2)$.  It is important to notice
that ${\bf Q}$ is not a multiplication operator in general.  Note that,
with suitable definition of domains, $B:\,Q_1 \to Q_2$ is equivalent
to $Q_1B^{*} = Q_1^{*}B^{*} = B^{*}Q_2^{*} = B^{*}Q_2,$ or $B^{*}:\,
Q_2 \to Q_1$.  Similarly ${\cal B} \sim B^{-1}$ satisfies $Q_1B^{-1} =
B^{-1}Q_2$, so ${\cal B}:\,Q_2 \to Q_1$.  We point out in passing
however that $B^{*}$ and $B^{-1}$ have opposite triangularities
(cf. (2.3), (2.5), etc.).  Now we will prove (note a priori $\Phi_1
\cap\Phi_2$ could be $\{0\}$ - cf. also Theorem 4.11)
\\[3mm]\indent {\bf THEOREM 3.1.}$\,\,$The operator V defined by
$$ {\cal Q}_2f(\lambda) = {\cal Q}_1(Vf)(\lambda)
\eqno(3.7)$$
\medskip
\noindent for $f \in A$ (below - $A$ dense in $D(Q_2)\cap H_1$) will
satisfy $VQ_2 f = Q_1 Vf$.
\\[2mm]\indent {\it Proof}:  Let $A = \{f \in H_1 \cap H_2;\,\,
(1 + |\lambda|){\cal Q}_2f(\lambda) \in L^2(d\Gamma_1) \cap
L^2(d\Gamma_2)\}$.  Note A will be dense under our assumptions since
the space of such $(1 + |\lambda|){\cal Q}_2f$ will be dense in
$L^2(d\Gamma_1)\cap L^2(d\Gamma_2)$.  For $f \in A$ one defines
$Vf$ via (3.7) and it follows via the diagram that
$$ \lambda{\cal Q}_2f(\lambda) = \lambda{\cal Q}_1(Vf)(\lambda)
\eqno(3.8)$$
\medskip
\noindent Now from $Q_i\psi_i = \lambda\psi_i$ one can say that for
$f \in \Phi_i \subset D(Q_i)$
$$ \lambda{\cal Q}_if(\lambda) = <f(x),\lambda\psi_i(x,\lambda)> =
\eqno(3.9)$$
$$= <f(x),Q_i\psi_i(x,\lambda)> = <Q_if,\psi_i(x,\lambda)> =
{\cal Q}_i(Q_if)$$
\medskip
\noindent and this extends to $f \in D(Q_i)$.  Now the left side of (3.8)
becomes (for $f \in A),\,\,{\cal Q}_2(Q_2f) =
{\cal Q}_1(VQ_2f)$ (by (3.7)) and the right side is ${\cal Q}_1(Q_1Vf)$,
provided $Vf \in D(Q_1)$.  But in fact, writing ${\cal Q}_i \sim
\bar{{\cal Q}}_i$ in $H_i$ via scalar products $(\,,\,)_i$, i.e.
$\bar{{\cal Q}}_i = (h,\psi_i)_i$, the equation ${\cal Q}_1(VQ_2 f) =
\lambda{\cal Q}_1(Vf)$ can be expressed as
$$ {\cal Q}_1(VQ_2f) = (VQ_2f,\psi_1(x,\lambda))_1 = (Vf,\lambda\psi_1)_1
= (Vf,Q_1\psi_1)_1
\eqno(3.10)$$
\medskip
\noindent (we emphasize $\psi_i \not\in H_i$ but in expressing the
action of $\bar{{\cal Q}}_i$ we will use the $(\,,\,)_i$ notation -
see section 4 for more detail).  Now elements $g \in D(Q_1)$ can be
expressed formally via $g = \int({\cal Q}_1g)\psi_1(x,\lambda)
d\Gamma_1$ so, with a little argument by approximation
$$ (VQ_2f,g)_1 = (VQ_2f,\int\hat{g}_1\psi_1d\Gamma_1)_1
\eqno(3.11)$$
$$= \int\hat{g}_1(VQ_2f,\psi_1)_1d\Gamma_1 = \int\hat{g}_1
(Vf,\lambda\psi_1)_1d\Gamma_1 = (Vf,Q_1g)_1$$
\medskip
\noindent This means $g \to (Vf,Q_1g)_1$ is continuous in the topology
of $H_1$ so $Vf \in D(Q_1^{*})$ and $VQ_2f = Q_1Vf\,\,(Q_1^{*} = Q_1)$.
A simpler argument can be based on $D(Q) = \{f:\,\hat{f}\,\,and\,\,
\lambda\hat{f}\in L^2(d\Gamma)\}$.
{\bf QED}
\\[3mm]\indent {\bf REMARK 3.3}$\,\,$ One could make various
assumptions regarding the $\Phi_i,\,\,D(Q_i),$ etc to produce a
cleaner looking theory.  If e.g. $H_1 = H_2 = H$ then
we have a traditional
transmutation framework.  We emphasize however that although
${\cal Q}_1:\,H \to L^2(d\Gamma_1)$ and ${\cal Q}_2:\,H \to
L^2(d\Gamma_2)$ are isometries, we cannot say that $V = {\cal Q}_1^{-1}
{\cal Q}_2:\, D(Q_2) \to D(Q_1)$ will extend to a bounded operator
in H.  This is simply because the injection $i: L = L^2(d\Gamma_1)\cap
L^2(d\Gamma_2) \to L^2(d\Gamma_1)$ may not be continuous when L has
the $L^2(d\Gamma_2)$ topology.  Thus in general V will not be a bounded
operator.  However ${\bf Q} = {\cal Q}_2{\cal Q}_1^{-1}:\, L^2(d\Gamma_1)
\to L^2(d\Gamma_2)$ will be continuous (recall the ${\cal Q}_i$ or
$\bar{{\cal Q}}_i$ are isometries).  Further if one has e.g.
$d\Gamma_1 = \gamma_1d\lambda,\,\,d\Gamma_2 = \gamma_2d\lambda,$
with $|\gamma_1/\gamma_2|\leq M<\infty$, then $i$ is bounded and the
theorem will apply for $f\in D(Q_2)$ ($V$ will be bounded in this
situation).

\section{Analysis of {\bf Q}}

We recall a theorem of Levitan [24] which states that every continuous
linear operator in a space of analytic functions ${\cal H}$ is locally
a linear differential operator of (possibly) infinite order.  Here
the appropriate topology is that of uniform convergence on compact
sets, i.e. $F_n \to F$ means that for any fixed compact K,
$sup_{\lambda \in K}|F_n(\lambda) - F(\lambda)| \to 0$ (we write this
as $F_n \stackrel{ucc}{\to} F$).  Note that a proof is easily constructed
via the Cauchy integral formula.  Thus let ${\bf Q}:\,{\cal H}(\Omega)
\stackrel{ucc}{\to}{\cal H}(\Omega)\,\,(\Omega$ open) be continuous
and linear and given a compact K with $z \in K$ let $K\subset   C\subset
\hat{K}\subset\Omega,\,\,\hat{K}$ compact (where C is a curve).  Then
${\bf Q}F_n(z) = (1/2\pi i)\oint_CF_n(\xi){\bf Q}(1/(\xi - z))d\xi$.
It follows that $F_n \to 0$ uniformly on $\hat{K}$ implies ${\bf Q}F_n
\to 0$ uniformly on K.  Consequently, locally
$$ {\bf Q}F(z) = \sum_0^{\infty}F^{(j)}(z)\frac{1}{2\pi i}\oint_C
(\zeta - z)^j{\bf Q}(1/(\zeta - z))d\zeta = \sum_0^{\infty}a_j(z)
\partial^jF(z)$$
\medskip
\noindent Note that ${\bf Q}(1/(\xi - z))$ must be defined here so
$1/(\xi -z)$ must be analytic for $\xi \in C,\,z \in K$.
\\[2mm]\indent Let now $W_i\subset H_i$ be the space of functions
such that ${\cal Q}_i(W_i)$ is entire when extended to ${\bf C}$.  One
thinks here of $C_0^{\infty}$ and Paley-Wiener (= PW) theorems for
example so many examples exist where $W_i$ will be dense (cf. [11;12]).
Let us also assume for convenience that $dM_{Q_i} = dx$ (so
$H_1 = H_2 = H)$.  A priori $\Phi_i$ and $W_i$ may not have any nice
relation but we note that the spaces $W^{\Omega}$ of [19] defined below
will usually be available as dense subspaces of $L^2(d\Gamma_1)\cap
L^2(d\Gamma_2)$ so that $\hat{\Phi}_i = {\cal Q}_i^{-1}W^{\Omega}
\subset H_i$ and $\hat{\Phi}'_i$ could be used for a Gelfand triple
(cf. remarks at the beginning of section 3).  Hence we will assume
$\Phi_i\subset W_i$ without loss of generality.  As for $W^{\Omega}$ we
define $W^{\Omega} = \{F\,\,entire;\,\,(1 + |\lambda|)^k|F(\lambda)|
\leq c_k exp(\Omega(b|Im(\lambda)|))\}$ where $\Omega(y) =
\int_0^y\xi(x)dx;\,\,\xi(x)\geq 0\,\,for\,\,\lambda > 0$.
$W^{\Omega}$ is a countably normed space with seminorms
$\|f\|_n = sup_{\lambda}(1 + |\lambda|)^n|F(\lambda)
exp(-\Omega(b|Im(\lambda)|))$ and the convergence of sequences is
defined by $F_n\stackrel{ucc}{\to}F$ with
$(1 + |\lambda|)^k|F_n(\lambda)| \leq c_kexp(\Omega(b|Im(\lambda)|))$
for all n,k.
In the particular case when $\xi(t) = 1$, i.e., $\Omega(y) = y$,
then $W^y = Z$, where Z is the space of entire functions of order one
and finite type (i.e. exponential type) defined by the family of
seminorms $\|f\|_n^Z = sup_{\lambda}|F(\lambda)|(1 + |\lambda|)^n
exp(-a|Im(\lambda)|)$.  Spaces of type W are known to be perfect and
in the analysis to follow we shall make occasional use of such spaces
(cf. [19] for more details).  An operator A is continuous in the space
$W^{\Omega}$ if it maps bounded sets into bounded sets or equivalently
if $f_n \to 0$ implies $Af_n \to 0$.  In the situations we consider
with Lebesgue-Stieljtes measures $d\Gamma$, if we assume $d\Gamma =
\gamma d\lambda;\,\,\gamma = O(|\lambda|^p)$, then $W^{\Omega}
\subset L^2(d\Gamma)$ since $F(\lambda) = O(|\lambda|^{-n})$ for
any $n > 0$ when $F \in W^{\Omega}$.  Now one has
\\[3mm]\indent{\bf THEOREM 4.1.}$\,\,$ Assume ${\bf Q}$ can be extended
to be a map ${\bf Q}:\,{\cal H}\to {\cal H}$, continuous in the ucc
topology and $H_1 = H_2 = H$.  If $f \in W_1$ then $\hat{f}_2(\lambda)
= {\bf Q}\hat{f}_1(\lambda)\in {\cal H}$ so $f\in W_2$ and locally
($D\sim d/d\lambda)$
$$ {\bf Q}\hat{f}_1(\lambda) = \hat{f}_2(\lambda) =
\sum_0^{\infty}a_n(\lambda)D^n\hat{f}_1(\lambda)
\eqno(4.1)$$
\medskip
\noindent It follows then, assuming $\Phi_i\subset W_i$ as discussed
above, that $\psi_i(x,\cdot)\in {\cal H}$ weakly and as transform objects
$\chi$ for ${\cal Q}$ acting via $f \to (f,\chi)$ for $f \in \Phi$
(or $f \in H$) one has (here $\partial^n_{\lambda}$ refers to a weak
or scalar derivative)
$$ \psi_2(x,\lambda) = \sum_0^{\infty}a_n(\lambda)\partial^n_{\lambda}
\psi_1(x,\lambda)
\eqno(4.2)$$
\\[2mm]\indent {\it Proof}:$\,\,$Here we say $\psi(x,\cdot)\in {\cal H}$
weakly if $\lambda \to <f(x),\psi(x,\lambda)> \in {\cal H}$ for any
$f\in \Phi$ (recall $\psi \in \Phi'$).  Now the formula (4.1) follows
from [23] since for $f\in\Phi_1\subset W_1$ we know $\hat{f}_1(\lambda)
= <f(x),\psi_1(x,\lambda)>\in {\cal H}$ (so $\psi_1\in {\cal H}$ weakly and
similarly $\psi_2\in {\cal H}$ weakly since $\Phi_2\subset W_2$).
The equation ${\bf Q}\hat{f}_1(\lambda) = \hat{f}_2(\lambda)$ can be
written formally as
$$ \hat{f}_2(\lambda) = \sum_0^{\infty}a_n(\lambda)\partial^n_{\lambda}
<f(x),\psi_1(x,\lambda)_1>
\eqno(4.3)$$
$$ = <f(x),\sum a_n(\lambda)\partial^n_{\lambda}
\psi_1(x,\lambda)>_1 = \int f(x)\psi_2(x,\lambda)dx$$
\medskip
\noindent We do not know if $f \in \Phi_2\,\,(f\in W_2)$ and even if
we assume $\Phi_1\cap \Phi_2$ is dense in $\Phi_2$, making the last
term in (4.3) $<f(x),\psi_2(x,\lambda)>_2$, this forces a comparison
of $<\,,\,>_1$ and $<\,,\,>_2$.  Hence we want to use $H = H_1 = H_2$
as an identification space and write integral signs in (4.3) instead
of $<\,,\,>_1$ (note however that we want to use $<\,,\,>_1$ first in
order to differentiate in $\lambda$ weakly).  This implies that as
transform objects $\chi$ acting via $f \to (f,\chi) = \Xi(f),\,\,f\in
\Phi$, we can make the identification (4.2).  ${\bf QED}$
\\[2mm]\indent There are a number of variations possible here (note also
Proposition 4.6 below which indicates that the continuity of ${\bf Q}$
is usually too strong).  We remark first however that by the Riesz
theorem one could also use $H'$ as an identification space in Theorem
4.4.  Thus $\psi\in\Phi'$ generates via ${\cal Q}$ a map $f\to
(f,\chi_{\psi}) = \Xi_{\psi}(f)$ where $\chi_{\psi}\sim\psi$ for
$f\in\Phi$ where $(f,\chi_{\psi}) = <f,\chi_{\psi}> = <f,\psi>$.  Since
a scalar product $\Xi_{\psi}(f) = (f,\chi_{\psi})\sim <f,\L_{\psi}>$
for $L_{\psi}\in H'$ we can identify $\psi\sim L_{\psi}\in H'$ etc.
in (4.2).  Another variation is to assume $\Phi = \Phi_1\cap\Phi_2$ is
dense in $\Phi_i$ with a suitable topology.  Then one can also use
$\Phi'$ as an identification space and write (4.2) as a genuine equation
in weak derivatives in $\Phi'$ (cf. Theorem 4.7).
Now regarding weak differentiability we recall that in the
dual of a barreled LCS (= locally convex topological vector space)
E the weak topology is equivalent to the topology of uniform convergence
on precompact sets in E.  If in addition bounded sets are relatively
compact in E (i.e. E is a Montel space) then the weak topology in $E'$
is equivalent to the strong topology.  Further the strong dual of a
Montel is Montel and evidently a barreled perfect space is Montel.
Noting that strict inductive limits of barreled spaces are barreled
one sees that Gelfand triples will often involve Montel spaces $\Phi$
and $\Phi'$.  Thus without great loss of generality we can assume
$\Phi_i\subset H\subset \Phi'_i$ is a Gelfand-Montel triplet.  Following
[16;30] we have then (note stronger theorems and a comprehensive
study of vector valued analytic functions are available in [21] but
we include Corollaries 4.2 and 4.3 for completeness and to facilitate
calculation)
\\[3mm]\indent {\bf COROLLARY 4.2}.$\,\,$ Let $\Phi_i\subset H\subset
\Phi'_i$ be Gelfand-Montel triples and assume the other hypotheses of
Theorem 4.1.  Then the $\psi_i$ are strongly analytic and the derivatives
$\partial^n_{\lambda}\psi_1$ in (4.2) represent strong derivatives.
\\[2mm]\indent {\it Proof:}$\,\,$ Write $\psi'_w$ for the weak derivative
$(\partial\psi/\partial\lambda)_w$ so that for any $f\in \Phi,\,\,
\partial_{\lambda}F(\lambda) = \partial_{\lambda}<f(x),\psi(x,\lambda)>
=<f(x),\psi'_w(x,\lambda)>$.  Let $B = \{[\psi(x,\lambda + \Delta\lambda)
-\psi(x,\lambda)]/\Delta\lambda\} = \{\Delta\psi/\Delta\lambda\}$ for
$|\Delta\lambda| < \epsilon$ say.  This is a weakly bounded set in $\Phi'$
since $\Delta\psi/\Delta\lambda - \psi'_w\to 0$ weakly.  Hence B is
bounded for the strong topology and evidently $\psi'_w\in \Phi =
\Phi''^{*} = \Phi^{*}$ (Montel spaces are reflexive).  Since $\psi'_w$
is weakly adherent to the bounded set $B\subset\Phi'$ we have
$\psi'_w\in\Phi'$.  Now use the fact that the weak and strong topologies
coincide on $B\cup \psi'_w$ to conclude that $\Delta\psi/\Delta\lambda
\to \psi'_w$ strongly.  ${\bf QED}$
\\[2mm]\indent Actually for analytic functions one has another
recourse based on the Cauchy integral formula.  Thus we know, since
$F(\lambda) = <f(x),\psi(x,\lambda)>\in {\cal H}$,
$$ <f(x),\psi(x,\lambda)> = \frac{1}{2\pi i}\oint_C\frac{<f(x),
\psi(x,\lambda)>}{\zeta - \lambda}d\zeta = <f(x),\frac{1}{2\pi i}
\oint_C\frac{\psi(x,\lambda)}{\zeta - \lambda}d\zeta
\eqno(4.4)$$
\medskip
\noindent (cf. [8;16] for vector valued integration - here C is e.g.
a circle in ${\bf C}$ around $\lambda$).  Assume $\Phi$ is reflexive
now in which case an integral such as $I(x,\lambda) = \frac{1}{2\pi i}
\oint_C[\psi(x,\lambda)/(\zeta - \lambda)]d\zeta$, with e.g.
$<f,\psi(x,\lambda)>$ continuous in $\lambda$, will belong to $\Phi''^{*}
=\Phi^{*}$.  If $\Phi'$ is complete or quasi-complete for example then in
fact $I(x,\lambda)\in \Phi'$.
Then for $f\in B$ bounded in $\Phi$ one obtains the equation
$$ <f,\frac{\Delta I}{\Delta\lambda} - I_{\lambda}> = <f,\frac{1}{2\pi i}
\oint_C\psi(x,\lambda)[\frac{1}{(\zeta - \lambda)(\zeta - \lambda
-\Delta\lambda)} - \frac{1}{(\zeta - \lambda)^2}]d\zeta>
\eqno(4.5)$$
\medskip
\noindent Now $\zeta\to\psi(x,\lambda)$ is weakly continuous so
$J = \{\psi(x,\zeta),\,\,|\zeta| = 1\}$ is weakly compact, hence
weakly bounded.  If $\Phi$ is barreled this means J is strongly
bounded and (4.5) shows $(\Delta I/\Delta\lambda) - I_{\lambda} \to 0$
strongly in $\Phi'$.  Hence, since reflexive implies barreled and
barreled spaces have quasicomplete duals (cf. [9]) we have
\\[3mm]\indent {\bf COROLLARY 4.3}$\,\,$ In addition to the hypotheses
of Theorem 4.1 assume the $\Phi_i$ are reflexive.  Then the $\psi_i$
are strongly analytic and the derivatives in (4.2) represent strong
derivatives.
\\[2mm]\indent Next we assume ${\bf Q}$ is continuous ${\cal H} \to
{\cal H}$ and $H_1 = H_2 = H$ with $\Phi_i\subset H$ being reflexive
and $\Phi_i\subset W_i$ so that Theorem 4.1 holds with Corollary 4.3.
Further by Theorem 3.1 for $f\in D(Q_2)$ one has $VQ_2f = Q_1Vf$.
A question now arises about the relation of Paley-Wiener properties
${\bf PWP}$ relative to $Q_1$ and $Q_2$.
We recall that in the model
situation of [11;12;28] for example the factorization ${\cal P}B^{*} =
{\cal Q}$ could be used. Thus
let us say that an operator P as above (modeled on second order
differential operators Q) has the ${\bf PWP}$ if (thinking of halfline problems
on $[0.\infty)$ for convenience of comparison to the examples of
section 2), ${\bf PWP}:\,\,$ $supp(f)\subset [0,\sigma]\Leftrightarrow
{\cal P}f = \hat{f}$ is even and entire in $k = \sqrt{\lambda}$
of exponential type $\sigma$ (i.e. $|\hat{f}(\lambda)| \leq cexp
(\sigma |Im(k)|)$).  Let ${\cal E}_{\sigma}$ be this space of $\hat{f}$.
In such an event, in the examples of section 2 the kernel $B:\,P\to Q$
is easily shown to be lower triangular ($Q\sim -D^2 + q$) so that
$B^{*}$ is given as in (2.3).  Then $supp(f)\subset [0,\sigma]
\Rightarrow supp(B^{*}f)\subset  [0,\sigma]$ and ${\cal P}B^{*}f
= {\cal Q}f\in {\cal E}_{\sigma}$.  Thus $\Rightarrow$ in ${\bf PWP}$
is transported from P to Q.  Let us now examine this in the present
context.  We remark that one could model our operators on differential
operators of order n by working with matrix differential operators but
there is also another recourse indicated in [5].  Thus one can often
rescale an operator Q by a suitable operator function T(Q) where
T(x) is to be defined by the requirements of the situation.  The spectral
measure $\Gamma(\lambda)$ becomes $\Gamma(T(\lambda))$ and
$\psi(x,\lambda)\to \psi(x,T^{-1}(\lambda))$.  This amounts to working
with $T^{-1}(\lambda)$ entire functions when discussing analyticity so
an entire $\hat{f}(\lambda) = <f(x),\psi(x,\lambda)> = \sum c_n\lambda^n
\to \sum c_n[T^{-1}(\lambda)]^{-n}$.  Thus for $T(x) = x^2,\,\,Cos(kx)
\to Cos(\sqrt{\lambda}x,\,\,k\in [0,\infty)\cup [0,i\infty)\to
\lambda\in (-\infty,\infty)$, etc.  One could then refer {\bf PWP} hypotheses
to some generic $\lambda$ etc. but we will keep the k framework here for
comparison to section 2.
\\[2mm]\indent Thus assume $Q_i$ as indicated and assume $Q_1$ has {\bf PWP}
(recall however that our technique recovers $V\sim B^{*}:\,Q_2 \to Q_1$ instead
of $B:\,Q_1 \to Q_2$ directly).  One obtains from (3.9) ${\cal Q}_2 f =
{\bf Q}{\cal Q}_1 f = {\cal Q}_1(Vf)$ and we assume ${\bf Q}:\,
{\cal E}_{\sigma}\to {\cal E}_{\sigma}$ algebraically.  Then ${\cal Q}_1 f$
even entire of exponential type $\sigma$ implies ${\cal Q}_2 f$ has
the same property along with ${\cal Q}_1(Vf)$.  Hence $supp(Vf)\subset
[0,\sigma]$ (which is consistent with our identification $V\sim B^{*}$
and (2.5)).  Let us (via experience in section 2) suppose V is an operator
$Vf(x) = \int_0^{\infty}V(x,y)f(y)dy$ with say V(x,y) continuous
in y for $y\ne x$ and deduce upper triangularity.  Thus (cf. [10])
let $supp(f)\subset [0,\sigma]\Rightarrow supp(\int_0^{\infty}V(x,y)
f(y)dy\subset [0,\sigma]$.  Let $f(y) = \delta_n(y-y_0)\in C_0^{\infty}
(0,\infty)\to\delta(y-y_0)$ so $Vf(x)\to  V(x,y_0)$ (assume such
$\delta_n$ belong to $D(Q_1)$ - via our discussion of $W^{\Omega}$
etc. this will pose no problem).  For $0 < y_0,\,\,supp(f)\subset
[0,y_0 +\epsilon]$ and we obtain $V(x,y_0) = 0$ if $x > y_0 + \epsilon$.
Thus $V(x,y) = 0$ for $x > y$ and
$$ Vf(x) = \int_x^{\infty}V(x,y)f(y)dy
\eqno(4.6)$$
\medskip
\noindent We have shown heuristically
\\[3mm]\indent {\bf PROPOSITION 4.4}.$\,\,$ Take $H_1 = H_2 = H = L^2(dx)$.
Assume $Q_1$ has {\bf PWP} and ${\bf Q}:\,{\cal E}_{\sigma}\to{\cal
E}_{\sigma}$
algebraically (no continuity is assumed) so $supp(f)\subset [0,\sigma]
\Rightarrow {\cal Q}_2f\in {\cal E}_{\sigma}$.  Let V defined by the diagram
after (3.6), $V = {\cal Q}_1^{-1}{\cal Q}_2$, be an operator $Vf(x) =
\int_0^{\infty}V(x,y)f(y)dy$ with kernel continuous in y for $y\ne x$ and
assume $D(Q_1)$ contains functions $\delta_n(y-y_0)\in C_0^{\infty}$.
Then V has the form (4.6).
\\[2mm]\indent We observe that other types of V are possible (e.g.
$f\to f$).  Now what about the converse, namely ${\cal Q}_2f\in
{\cal E}_{\sigma}\Rightarrow supp(f)\subset [0,\sigma]$?  Naively
one expects an operator $V^{-1}$ to be upper triangular of the same
type as V, given that $V^{-1}$ exists.  This seems to be a difficult
question however and the constructions in [11;12;23;25;26;28] are
based on completeness theorems etc. for the generalized eigenfunctions.
We will turn to this matter now.  The condition ${\bf WLP}$ (weak local
property): $supp(f)\subset [0,\sigma]\Leftrightarrow supp(Vf)\subset
[0,\sigma]$
can be used to describe this situation.  Thus {\bf WLP} $\sim$ V upper
triangular $\Leftrightarrow V^{-1}$ upper triangular, but the idea is
more general and one can shortcut the development of eigenfunction
machinery and achieve on the face of it greater generality by using WLP
as a hypothesis in various theorems.  For example from what we have
stated there follows.
\\[3mm]\indent {\bf THEOREM 4.5}:$\,\,$ Assume $Q_1$ has ${\bf PWP}$
and ${\bf Q}:\,{\cal E}_{\sigma}\to {\cal E}_{\sigma}$ algebraically
(no continuity involved) with $H_1 = H_2 = H = L^2$.  Then $supp(f)
\subset [0,\sigma]\Rightarrow supp(Vf)\subset [0,\sigma]$ and
${\cal Q}_2f\in {\cal E}_{\sigma}$.  Further if V has ${\bf WLP}$ then
${\cal Q}_2f\in {\cal E}_{\sigma}\Rightarrow supp(Vf)\subset [0,\sigma]
\Rightarrow supp(f)\subset [0,\sigma]$ so $Q_2$ has ${\bf PWP}$.
On the other hand suppose $Q_1$ and $Q_2$ have ${\bf PWP}$.  Then
immediately from the diagram after (3.6),
$supp(f)\subset [0,\sigma] \Rightarrow supp(Vf)\subset [0,\sigma]$ and
$supp(Vf)\subset [0,\sigma] \Rightarrow supp(f)\subset [0,\sigma]$
so V has ${\bf WLP}$.
\\[2mm]\indent Thus let us look now at ${\cal Q}_2f = {\cal Q}_1(Vf)\,\,
(H_1 = H_2 = H = L^2)$ and try to express the {\bf WLP} in terms of
generalized eigenfunctions.  Written out this says formally, for
suitable f (cf. (4.3))
$$ \int f(x)\psi_2(x,\lambda)dx = \int (Vf)(x)\psi_1(x,\lambda)dx =
\int f(x)V^{*}\psi_1(x,\lambda)dx
\eqno(4.7)$$
\medskip
\noindent Here one imagines e.g. V as an integral operator such as
(4.6) (or better (2.3)) with $V^{*}$ the $L^2$ adjoint which is
presumed to be able to act on $\psi_1\in\Phi'_1$.  Now one knows that in
examples from [11;12;23;25;26;28] based on Cos(kx) etc. the generalized
eigenfunctions $\psi_i(x,\lambda)$ will be often $C^0,\,C^2,$ or
$C^{\infty}$ in x and analytic in $\lambda$ (so the hypotheses on
{\bf Q} in Theorem 4.1 are not a priori unreasonable).  Moreover
$V^{*}\sim B^{**} = B$ should have the form (via (2.3))
$$ V^{*}f(y) = f(y) + \int_0^yV(x,y)f(x)dx
\eqno(4.8)$$
\medskip
\noindent so it is not unrealistic to expect (4.7) with some interpretation
to yield
$$ \psi_2(x,\lambda) = (V^{*}\psi_1)(x,\lambda) = \psi_1(x,\lambda)
+ \int_0^xV(y,x)\psi_1(y,\lambda)dy
\eqno(4.9)$$
\medskip
\noindent This is the natural context then and we want to see how
much can be deduced in the more abstract situations.
\\[2mm]\indent From [19] (cf. also [4;18]) we have formally for
$h\in H = L^2$
$$ h(x) = \int(\psi(y,\lambda),h(y))\psi(x,\lambda)d\Gamma;
\eqno(4.10)$$
$$\|h\|^2 =\int |(\psi,h)|^2d\Gamma; (\psi(y,\lambda),h(y)) = 0\,\,
\forall\lambda\Rightarrow h = 0$$
\medskip
\noindent From this one writes formally
$$ \delta(x-y) = \int\psi(x,\lambda)\psi(y,\lambda)d\Gamma
\eqno(4.11)$$
\medskip
\noindent whose meaning is specified by (4.10).  Further we can formally
represent kernels as in section 2 via
$$ \beta(y,x) = kerB = \int\psi_1(x,\lambda)\psi_2(y,\lambda)d\Gamma_1;
\eqno(4.12)$$
$$ \gamma(x,y) = ker{\cal B} = \int\psi_1(x,\lambda)\psi_2(y,\lambda)
d\Gamma_2$$
\medskip
\noindent In classical situations these are integrals such as
$\frac{2}{\pi}\int_0^{\infty}Cos(\lambda x)Cos(\lambda y)d\lambda =
\delta_{+}(x-y)$ (half line delta function - cf. [11;12]) which are not
strictly well defined integrals but acquire a meaning via distribution
theory.  For singular differential operators one can find many formulas
of the form (4.12) in [11;12] with standard special functions for the
$\psi_i$ where everything makes good sense via distribution theory.
If formulas such as (4.10) for example cause anxiety one can think of
approximating $h\in H$ by $\phi_n\in \Phi$ and/or realize that
$h(x) = \int(\psi(y,\lambda),h(y))\psi(x,\lambda)d\Gamma$ is simply
another way of saying $h = {\cal Q}^{-1}{\cal Q}h$.  Similarly (4.12)
says e.g. (cf. (2.8), (2.10))
$$ Bf(y) = <\beta(y,x),f(x)> = \int{\cal Q}_1f\psi_2(y,\lambda)
d\Gamma_1 = \tilde{{\cal Q}}_2{\cal Q}_1f(y);$$
$$ \tilde{{\cal Q}}_2F(y) = \int F(\lambda)\psi_2(y,\lambda)d\Gamma_1
\eqno(4.13)$$
\medskip
\noindent We see formally that $\gamma$ in (4.12) represents an inverse
kernel as follows.  First for simplicity assume the $d\Gamma_i$ are
absolutely continuous with $d\Gamma_i = \gamma_i(\lambda)d\lambda$.  Then
to go with (4.11) one should have formally
$$ \int \psi(x,\lambda)\psi(x,\mu)dx = \delta(\lambda - \mu)/\gamma(\mu)
\eqno(4.14)$$
\medskip
\noindent This is equivalent formally to
$$ F(\lambda) = \int[\delta(\lambda - \mu)/\gamma(\mu)]F(\mu)d\Gamma(\mu) =
\eqno(4.15)$$
$$ = \int\psi(x,\lambda)(\int\psi(x,\mu)F(\mu)d\Gamma(\mu))dx =
{\cal Q}{\cal Q}^{-1} F(\lambda)$$
\medskip
\noindent just as (4.11) is formally equivalent to
$$ f(x) = \int\delta(x-y)f(y)dy =
\eqno(4.16)$$
$$ = \int\psi(x,\lambda)(\int\psi(y,\lambda)f(y)dy)d\Gamma(\lambda)
= {\cal Q}^{-1}{\cal Q}f$$
\medskip
\noindent Given the formal structure indicated one has e.g.
$$ <\gamma(x,y),\beta(y,\xi)> = \int\psi_1(x,\lambda)(\int\psi_2
(y,\mu)dy)\psi_1(\xi,\mu)d\Gamma_2(\lambda)d\Gamma_1(\mu) =
\eqno(4.17)$$
$$ = \int\psi_1(x,\lambda)\frac{\delta(\lambda - \mu)}{\gamma_2(\mu)}
\psi_1(\xi,\mu)\gamma_2(\lambda)d\lambda\gamma_1(\mu)d\mu =
\int\psi_1(x,\lambda)\psi_1(\xi,\lambda)\gamma_1(\lambda)d\lambda =
\delta(x - \xi)$$
\medskip
\noindent so the inversion kernels as in (4.12) are natural.
\\[2mm]\indent Now triangularity in the classical examples is
proved via hyperbolic PDE or via analyticity properties of the
generalized eigenfunctions (cf. [11;12]).  Further the lower
triangularity of B and ${\cal B}$ is equivalent to the {\bf WLP}
for $B^{*}\sim V$ for the classical examples.  However both the
contour integration technique from [11;12] or PDE techniques as
in [11;12;23;25;26;28] require more detailed knowledge of either
the $\psi_i$ or the $Q_i$.  Consider now as a prototypical Q the
operator $Q_1 = -D^2$ in $L^2(0,\infty)$ with the generalized
eigenfunctions $\psi_1(x,\lambda) = Cos(kx)\,\,(\lambda = k^2)$.  Let
$Q_2$ have {\bf PWP} (as is common in examples) and look again at
Theorem 4.1 (plus the corollaries).  Suppose ${\bf Q}:\,{\cal H}\to
{\cal H}$ is continuous in the ucc topology, leading to (4.2).  This
would imply formally $(\partial_{\lambda} = \frac{1}{2k}\partial_k)$
$$ \psi_2(x,\lambda) \sim \sum_0^{\infty}a_n(k^2)(\frac{\partial_k}{2k})^n
Cos(kx) =
\eqno(4.18)$$
$$ = \sum_0^{\infty}a_n(k^2)[P_n(x,\frac{1}{k})Cos(kx) + \hat{P}_n
(x,\frac{1}{k})Sin(kx)] $$
\medskip
\noindent Thus formally at least $\psi_2$ would be analytic in x which
is unlikely with examples like $Q_2 = -D^2 + q(x)$ with q only
continuous.  Hence {\bf Q} cannot be continuous as indicated in general
and we emphasize this via
\\[3mm]\indent {\bf PROPOSITION 4.6}.$\,\,$ The hypothesis ${\bf Q}:
{\cal H}\to {\cal H}$ continuous in the ucc topology in Theorem 4.1 is
usually too strong.
\\[2mm]\indent Let us see if $Q_1$ and $Q_2$ close in some sense will
imply {\bf Q} continuous in some sense (perhaps not ${\cal H}\to
{\cal H}$ in ucc topology but in a topology one can adapt to the
Levitan theorem).  Note that even though $Q\psi = \lambda\psi,\,\,
\psi\in\Phi'$, this only holds for $\lambda\in\sigma_Q\subset {\bf R}$
and one cannot directly generate an analyticity argument in $\lambda$.
Hence we will have to assume $\psi_i(x,\lambda)\in {\cal H}$ in
$\lambda$ for the next result.  First in order to compare $Q_1$ and
$Q_2$ let assume $H_1 = H_2 = H$ and $\Phi = \Phi_1\cap\Phi_2$ is
dense in $\Phi_i$.  Put on $\Phi$ the topology of simultaneous convergence
in $\Phi_1$ and $\Phi_2$ so $\Phi\subset\Phi_i\subset H\subset\Phi'_i
\subset\Phi'$.  Then $\Phi'$ can be used as an identification space
for the $\psi_i$.  Also $\lambda\to\psi_i(x,\lambda)$ entire with values
in $\Phi'_i$ will imply these functions are entire with values in
$\Phi'$.  Further without great loss of generality one could assume
also e.g. that ${\cal Q}_i^{-1}Z = \Phi_i$ (see remarks before Theorem 4.1).
Next we consider the possibility of ${\bf Q}: Z\to {\cal H}$ being
continuous (the reason for this will emerge in the proof of Theorem
4.7 to follow).  Recall that $F_n\to 0$ in $Z\sim \|F_n\|_m^Z\to 0$
where $\|F\|_m^Z = sup|F|(1 + |\lambda|)^m exp(-a|Im\lambda|)$
(adjust here $\lambda\sim k^2$ as needed for Q modeled on $D^2$ etc. -
cf. remarks before Proposition 4.4).  Now locally we can still write
${\bf Q}F(z) = (1/2\pi i)\oint_C F(\xi){\bf Q}(\frac{1}{\xi - z})d\xi$
for $z\in K\subset C\subset \hat{K}\subset \Omega$ as before and
$F_n\to 0$ in Z implies $F_n\to 0$ in ${\cal H}$ so our previous
discussion applies.  The only additional feature is that $sup|F_n|
(1 + |\lambda|)^m\leq k_m$ for $\lambda \in {\bf R}$ and any m, so if
$d\Gamma_1 = \gamma_1d\lambda$ with $\gamma_1 = O(|\lambda|^p)$
for some p, then
$$ \|F_n\|^2_{L^2(d\Gamma_1)} = \int|F_n|^2\frac{|\lambda^p}
{1 + |\lambda|)^{2m}}(1 + |\lambda|)^{2m}d\lambda \leq
\eqno(4.19)$$
$$ \int\frac{|\lambda|^p}{1 + |\lambda|)^{2m}}(|F_n|(1 +
|\lambda|)^m)^2d\lambda
\leq k_m^2\|F_n\|_m^Z\,\,\,(2m \geq p+2)$$
\\[3mm]\indent {\bf THEOREM 4.7}$\,\,$ Assume $H_1 = H_2 = H =
L^2(dx)$ for convenience, $\Phi = {\cal Q}_1^{-1}Z\cap {\cal Q}_2^{-1}Z
= \Phi_1\cap\Phi_2$ is dense in $\Phi_i$, and in H, and that
$\psi_i(x,\lambda)\in\Phi'$ are entire in $\lambda$.  Set $r(x,\lambda) =
\psi_2(x,\lambda) - \psi_1(x,\lambda)$ in $\Phi'$ with $\hat{r}_1(\nu,
\lambda) = \int r(x,\lambda)\psi_1(x,\nu)dx$.  Assume (A) For K compact
there exists $g(\nu,K) \geq 0\in L^2(d\Gamma_1(\nu))$ such that $|\hat{r}_1
(\nu,\lambda)| \leq g(\nu,K)$ for $\lambda\in K$ (B) $d\Gamma_1$ is
absolutely continuous with $d\Gamma_1 = \gamma_1d\lambda$ and
$|\gamma_1| = O(|\lambda|^p)$ (C) ${\bf Q}: Z\to {\cal H}$ is
continuous.  Then (4.2) holds locally in $\Phi'$.
\\[2mm]\indent {\it Proof}:$\,\,$ Clearly $r(x,\lambda)$ is entire
in $\lambda$ with values in $\Phi'$ and via Parseval for ${\cal Q}_1$ one has
$$ {\bf Q}\hat{f}_1(\lambda) = \hat{f}_2(\lambda) = \hat{f}_1(\lambda) +
\int\hat{r}_1(\nu,\lambda)\hat{f}_1(\nu)d\Gamma_1(\nu)
\eqno(4.20)$$
\medskip
\noindent (i.e. $\int f(x)r(x,\lambda)dx = \int \hat{f}_1\hat{r}_1
(\nu,\lambda)d\Gamma_1(\nu))$.  This is an integral equation of
Carleman type and (A) ensures that one can use dominated convergence
ideas in the integration.  In particular $\hat{f}_2(\lambda)$ and
$\hat{f}_1(\lambda)$ are entire and we will show that $\hat{f}_n^1
\stackrel{Z}{\to} 0 \Rightarrow {\bf Q}\hat{f}_n^1 = \hat{f}_n^2
\stackrel{ucc}{\to} 0$.  One need only estimate the integral term
$J(\lambda) = \int \hat{r}_1(\nu,\lambda)\hat{f}_1(\nu)d\Gamma_1(\nu)$.
First for $f_n\in\Phi$ as above, one has $|\hat{f}_n^1(\nu)|\leq k_m
(1 + |\nu|)^{-m}$ for any m as above and we choose m so that
$\|\hat{f}_n^1\|_{L^2(d\Gamma_1}\leq k_m$.  Then consider for $\hat{f}_n^1
\in Z$
$$ sup_{\lambda\in K}|J_n(\lambda)| \leq \int g(\nu,K)|\hat{f}_n^1(\nu)|
d\Gamma_1(\nu)
\eqno(4.21)$$
$$ \leq (\int g^2d\Gamma_1)^{\frac{1}{2}}(\int |\hat{f}_n^1(\nu)|^2
d\Gamma_1)^{\frac{1}{2}} \leq Gk_m\|\hat{f}_n^1\|_m\,\,(2m \geq p+2)$$
\medskip
\noindent But $\hat{f}_n^1\stackrel{Z}{\to} 0$ means $\|\hat{f}_n^1\|_m
\to 0$ for any m, hence in particular for $2m\geq p+2,\,\,J_n(\lambda)
\to 0$ in ${\cal H}$.  Consequently ${\bf Q}\hat{f}_1 = \hat{f}_2$ can
be written as in (4.1) for $f\in\Phi$.  This leads to (4.3) and the
identification as in (4.2) with the $\psi_i$ as transform objects,
or locally in $\Phi'$ strongly.  {\bf QED}
\\[3mm]\indent {\bf REMARK 4.8}:$\,\,$ The proof holds for any
$W^{\Omega}$ type space with adjustment of hypotheses on $\Phi_i$
(i.e. $\Phi_i = {\cal Q}_i^{-1}W^{\Omega}$, etc.) but the strongest
version would involve only ${\bf Q}:\,{\cal H}\to {\cal H}$ algebraically
and ${\bf Q}F_n\to 0$ in ${\cal H}$ provided $F_n\to 0$ in ${\cal H}$
plus $(1 + |\lambda|)^m sup_{\lambda\in {\bf R}}|F_n(\lambda)|\to 0$
for some $2m\geq p+2$.  Of course other hypotheses on $|\hat{r}_1(\nu,
\lambda)|$ could also be made.  Note that given $\Phi$ as indicated
in Theorem 4.7, (4.18) would suggest that for $\psi_1\sim Cos(kx)$,
which is analytic in x, if $Q_2$ is suffficiently close to $Q_1 = -D^2$
(measured by $\hat{r}_1$), then $\psi_2(x,\lambda)$ might also be
analytic in x, via continuity of {\bf Q} in some sense.  Further
investigation of such matters is clearly indicated.
\\[3mm]\indent {\bf REMARK 4.9}.$\,\,$ Given $Q_1 = -D^2$ and $Q_2 =
-D^2 + q(x)$ one has formally $-\psi''_2 + q\psi_2 = \lambda\psi_2 =
k^2\psi_2$.  Given a transmutation B as in (2.2) with $\psi_2(x,\lambda)
=Cos(kx) + \int_0^x K(x,t)Cos(kt)dt$ one will have $\psi_2$ analytic
in $\lambda$ or k.  Further from the differential equation, if q(x)
is analytic in x one expects $\psi_2$ will be analytic in x.  Thus such
a situation should involve formulas of the type (4.18).  In other words
analytic q(x) could produce closeness of $Q_2$ and $Q_1$ in the sense
of Theorem 4.7.  One suspects also that relations could be established
here to the results of [27].
\\[3mm]\indent {\bf REMARK 4.10}.$\,\,$ We note here another approach
to determining transmutations analogous to the diagram after (3.6).
Thus ${\cal B}:\,Q_2\to Q_1$ can be formally defined as
${\cal B} = \tilde{{\cal Q}}_1{\cal Q}_2$ (cf. (2.11), (4.13)), where
$\tilde{{\cal Q}}_1F(x) = \int F(\lambda)\psi_1(x,\lambda)d\Gamma_2$.
Consider for simplicity a situation where say $Z\subset L^2(d\Gamma_1)
\cap L^2(d\Gamma_2)$ is dense in each $L^2$ and $\Phi_i = {\cal Q}_i^{-1}Z$
with $\Phi = \Phi_1\cap \Phi_2$ dense in $H = H_1 = H_2$
and in $D(Q_i)$.  Then look
at a diagram ($d\Gamma_i = \gamma_i d\lambda$) with ${\bf P} =
\gamma_2/\gamma_1$
\begin{eqnarray}
\setlength{\jot}{6pt}
\hat{f}_2 & \stackrel{{\bf P}}{\to} & {\bf P}\hat{f}_2\nonumber\\
{\cal Q}_2\uparrow & \searrow \tilde{{\cal Q}}_1 & \downarrow {\cal Q}_1^{-1}
\nonumber\\
f\in \Phi & \stackrel{{\cal B}}{\to} & {\cal B}f\in H\nonumber
\end{eqnarray}
\medskip
\noindent and observe here that for ${\bf P} = \gamma_2/\gamma_1$
$$ \int\hat{f}_2\psi_1 d\Gamma_2 = \int\hat{f}_2\frac{\gamma_2}{\gamma_1}
\psi_1 d\Gamma_1; \tilde{{\cal Q}}_1 = {\cal Q}_1^{-1}{\bf P};
{\cal B} = \tilde{{\cal Q}}_1{\cal Q}_2
\eqno(4.22)$$
\medskip
\noindent We record here in this connection
\\[3mm]\indent {\bf THEOREM 4.11}.$\,\,$ For ${\bf P} = \sqrt{\gamma_2/
\gamma_1}$ the diagram (with $\tilde{{\cal Q}}_1$ removed)
defines a new kind of transmutation
${\cal B}:\, Q_2\to Q_1$, satisfying $Q_1{\cal B}f = {\cal B}Q_2 f$
for $f\in D(Q_2)$.  For ${\bf P} = \gamma_2/\gamma_1$ with $|\gamma_2/
\gamma_1|\leq c$ (and $\tilde{{\cal Q}}_1$)
we get ${\cal B}\sim B^{-1}$ (cf. (2.11) and (4.22)).
\\[2mm]\indent {\it Proof}:$\,\,$ Take $f\in D(Q_2)$ so that
$$ \lambda {\cal Q}_1({\cal B}f) = {\bf P}\lambda {\cal Q}_2(f) =
{\bf P} {\cal Q}_2(Q_2f) = {\cal Q}_1({\cal B}Q_2 f)
\eqno(4.23)$$
\medskip
\noindent (cf. here (3.8)-(3.11) - the conditions on ${\bf P}$ insure
that everything makes sense).
Then we want to show $({\cal B}Q_2f,h)
=({\cal B}f,Q_1h)$ for $h\in D(Q_1)$ which would imply that
${\cal B}f\in D(Q_1)$ with ${\cal B}Q_2 f = Q_1 {\cal B}f$.  Thus
represent h as $h = {\cal Q}_1^{-1}{\cal Q}_1h = \int ({\cal Q}_1h)
\psi_1\gamma_1 d\lambda$ and one has $Q_1h = \int ({\cal Q}_1 h)\lambda
\psi_1\gamma_1 d\lambda$ with
$$ (Q_1 h,{\cal B}f) = \int ({\cal Q}_1 h)\gamma_1\lambda (\psi_1,
{\cal B}f) d\lambda =
\eqno(4.24)$$
$$ = \int ({\cal Q}_1 h)\gamma_1\lambda {\cal Q}_1 ({\cal B} f)d\lambda
=\int ({\cal Q}_1 h)\gamma_1 {\cal Q}_1({\cal B}Q_2 f)d\lambda =
(h,{\cal B}Q_2f)$$
\medskip
\noindent (the last equation by Parseval).  {\bf QED}
\\[2mm]
\noindent This shows that ${\cal B}:\,Q_2\to Q_1$ is a transmutation and
we have developed constructions of ${\cal B}\sim B^{-1}$ and $V\sim B^{*}$.
Note also that if {\bf P} is a polynomial for example then
{\bf P} maps ${\cal H}\to {\cal H}$ or $Z\to Z$.  One can surely enhance
the investigation of {\bf PWP} and
related matters using Theorem 4.11 with the previous results
and we will return to this at another time.

\end{document}